\documentclass{llncs}
\usepackage{amsmath}
\usepackage{tabularx}
\usepackage{graphicx,rotating}
\usepackage{tikz}
\usetikzlibrary{matrix}

\mathchardef\mhyphen="2D 

\begin{document}

\title{ Set Matching: An Enhancement of the Hales-Jewett Pairing Strategy}
\author{Jos W.H.M. Uiterwijk}
\institute{Department of Data Science and Knowledge Engineering \\ Maastricht University \\
\email{uiterwijk@maastrichtuniversity.nl}}
\date{}

\maketitle

\begin{abstract}
When solving $k$-in-a-Row games,  the Hales-Jewett pairing strategy \cite{HalJew1963} is a well-known strategy to prove that specific positions are (at most) a draw. It requires two empty squares per possible winning line (group) to be marked, i.e., with a coverage ratio of 2.0.

In this paper we present a new strategy, called {\em Set Matching}. A matching set consists of a set of nodes (the markers), a set of possible winning lines (the groups), and a coverage set indicating how all groups are covered after every first initial move. This strategy needs less than two markers per group. As such it is able to prove positions in $k$-in-a-Row games to be draws, which cannot be proven using the Hales-Jewett pairing strategy.

We show several efficient configurations with their matching sets. These include Cycle Configurations, BiCycle Configurations, and PolyCycle Configurations involving more than two cycles. Depending on configuration, the coverage ratio can be reduced to 1.14.

Many examples in the domain of solving $k$-in-a-Row games are given, including the direct proof (without further investigation) that the empty $4 \times 4$ board is a draw for 4-in-a-Row.
\end{abstract}

\section{Introduction}

Playing $k$-in-a-Row games is a popular pastime activity \cite{Gardner1983}. They are played by two players alternately claiming a square, trying to form $k$   squares of their color in a straight line on some rectangular $m \times n$ board. Therefore, these games are also denoted as $mnk$-games \cite{UitHer2000}. As such $k$-in-a-Row games are a subset of {\em strong positional games} \cite{Beck2005,Beck2008}. A useful tool in solving $mnk$-games is the Hales-Jewett (HJ)-pairing strategy \cite{HalJew1963}, by which certain positions can be proven to be at most a draw for the first player.   The disadvantage is that for every group involved two distinct empty squares  should be used. 

In this paper we present results of a study to transform the pairing strategy into a method where a set of groups can be proven to be at most a draw using less than two empty squares  per group.  We denote such a set as a {\em matching set}. In Section 2 we give some background theory. Then three classes of matching sets using one, two, and more cycles are presented in Sections 3-5. The configurations are exemplified by positions from 4-in-a-Row games. In Section 6 we provide an overview of the configurations and their efficiencies, and give some experimental results. Conclusions and an outlook to future research are given in Section 7.

\section{Background}

We first give a theoretical framework for $mnk$-games and provide some useful notions. We then explain the Hales-Jewett pairing strategy and introduce our new method of Set Matching.

\begin{definition}
{\em (Taken from  \cite{Beck2005,Beck2008})} A {\em  positional game} is a hypergraph ($X$, $H$), where the set $X$  contains nodes forming the game board  and $H  \subseteq 2^X$ is a family of target subsets of $X$. During the game, two players alternately select one previously unclaimed element of the board. When the goal of the game for both players is to be the first to  claim all elements of a target subset, the game is  a {\em strong positional game}.
\end{definition}
 
\noindent $mnk$-Games are examples of strong positional games. The target subsets are all the possible winning lines, also called {\em groups}, of $k$ squares in a straight line. The first player who claims all elements of a group wins the game. If no player achieves a win, the game is a draw.
Groups can interact in several ways.
\begin{definition}
When groups have common nodes, they are called {\em intersecting}. The common nodes are denoted as {\em intersections} or {\em corner nodes}. Nodes in a group not being intersections are called {\em edge nodes}.
A pair of groups with more than one intersection is called {\em overlapping}.
\end{definition}

In the remainder of this paper we only consider sets of intersecting but not overlapping groups. 
A useful theorem for strong positional games, due to John Nash, but probably first published in \cite{HalJew1963}, is known as the {\em strategy-stealing argument}, which results in an $mnk$-game being either a first-player  win or a draw.
 To prove that a set of non-overlapping groups is at most a draw for the first player, the well-known HJ-pairing strategy can be used.
\begin{definition}
A {\em Hales-Jewett  (HJ)-pairing} for a set of groups is an assignment of pairs of empty nodes  ({\em markers}) to all groups, such that every group is covered by a marker pair. The second player hereby guarantees for every marker of a pair played by the first player to respond immediately by playing the second marker of the pair, thus covering every group of the set by at least one second-player stone, preventing  a possible win by the first player in that group. 
\end{definition}

As an enhancement of a HJ-pairing we introduce matching sets, needing less markers than 2 per group to cover the whole set of groups. A formal notation for a matching set $M$ is as follows. 
\begin{definition}
A {\em matching set} $M$ is a triple $(N, G, C)$, where $N$ denotes the set of marker nodes used for the matching, $G$ the set of groups covered by the matching, and $C$ the  set of coverings. A covering consists of a black move followed ($\rightarrow$) by the white response, followed by a matching set for the remaining groups not covered by the white response. 
\end{definition}

We use the following conventions. 
1) Common abbreviations will be used, like $a$-$h$ for the markers $a$ to $h$. 
2) Every group is denoted by all the markers in the group, where for brevity reasons we  omit parentheses and commas. For instance, the group containing markers $a$ and $b$ is denoted by just $ab$. 
3) We omit coverings that are equivalent to other mentioned coverings (mostly symmetric in the configuration sketch). 
4) In case a matching set just consists of a naive HJ-pairing, we just give a set with marker pairs.

\section{Set Matching using Cycle Configurations}

In this section we give several configurations that need less than 2 markers per group. For each such configuration we give a sketch and an actual board example taken from a 4-in-a-Row game. We always assume that the positions are with Black to move and that White shows that Black cannot win the position. The example positions are chosen such that no other groups are present except those in the matching set, that all empty squares in the  groups that are not used as markers are Black (to prevent simple HJ-pairings), that there are as many White as Black stones (i.e., legal positions with Black to move), and that there are no forced lines outside the groups.

\subsection{The Triangle Configuration}

\begin{figure}[!h]
\begin{center}
\begin{tikzpicture}[scale=0.7]
\draw[fill=black] (-1.0,-1.732) circle (3.0pt);
\node at (-1.3,-1.732) {$a$};
\draw[fill=black] (-1.0,1.732) circle (3.0pt);
\node at (-1.3,1.732) {$b$};
\draw[fill=black] (2.0,0.0) circle (3.0pt);
\node at (2.3,0.0) {$c$};
\draw[fill=black] (0.5,-0.866) circle (3.0pt);
\node at (0.5,-1.266) {$d$};
\draw[fill=black] (0.5,0.866) circle (3.0pt);
\node at (0.5,1.266) {$e$};
\draw[dashed] (2.0,0.0) -- (-1.0,1.732);
\draw[dashed] (2.0,0.0) -- (-1.0,-1.732);
\draw[dashed] (-1.0,1.732) -- (-1.0,-1.732);
\end{tikzpicture} 
\hspace*{1.0cm}
\begin{tikzpicture}[scale=0.9]
\draw[fill=black] (-1.4,-1.05) circle (2.0pt);
\draw[fill=black] (-1.4,1.05) circle (2.0pt);
\draw[fill=black] (0.0,-1.05) circle (2.0pt);
\draw[fill=black] (0.0,-0.35) circle (2.0pt);
\draw[fill=black] (0.7,-1.05) circle (2.0pt);
\draw[dashed] (-1.4,-1.05) -- (-1.4,1.05);
\draw[dashed] (-1.4,1.05) -- (0.7,-1.05);
\draw[dashed] (-1.4,-1.05) -- (0.7,-1.05);
\input{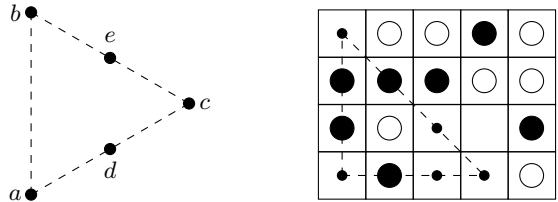}
\matrix[square matrix=.63cm]
{
	&	& 	& 	& 	\\
	& 	&	& 	& 	\\
	& 	& 	&	& 	\\
	& 	& 	& 	&	\\
};
\draw[fill=black] (-0.7,-1.05) circle (5pt);
\draw (1.4,-1.05) circle (5pt);
\draw[fill=black] (-1.4,-0.35) circle (5pt);
\draw (-0.7,-0.35) circle (5pt);
\draw[fill=black] (1.4,-0.35) circle (5pt);
\draw[fill=black] (-1.4,0.35) circle (5pt);
\draw[fill=black] (-0.7,0.35) circle (5pt);
\draw[fill=black] (0.0,0.35) circle (5pt);
\draw (0.7,0.35) circle (5pt);
\draw (1.4,0.35) circle (5pt);
\draw (-0.7,1.05) circle (5pt);
\draw (0.0,1.05) circle (5pt);
\draw[fill=black] (0.7,1.05) circle (5pt);
\draw (1.4,1.05) circle (5pt);
\end{tikzpicture}
\end{center}
\vspace*{-0.2cm}
\caption{Sketch and example position  illustrating the Triangle Configuration.}
    \label{fig:Triangle}
\end{figure}

To achieve a reduction of  the number of markers needed, we first investigated the idea that the groups intersect in a cyclic way. The smallest such set of groups is when three such groups are involved.
The Triangle Configuration  is formed by three pair-wise intersecting groups. All three intersections are markers. Further two of the three groups have one additional marker. See Fig. \ref{fig:Triangle} (left) for a sketch. Markers are required to be empty. Non-marker nodes can be  empty or black. 

\vspace*{-0.3cm}

In this configuration White always responds with $a$ or $b$ after any black marker move. We denote $a$ and $b$ as the {\em main markers}. By doing this, two of the three groups will be covered and the third group then can be covered by  two unused markers in that group using a naive HJ-pairing.  The formal notation for a matching set (not necessarily unique) of this configuration is:
\begin{align*}
M = ( \{&a \mhyphen e\}, \{ab, adc, bec\}, \\
\{ &a \rightarrow b \{(c,d)\}, c \rightarrow a \{(b,e)\}, d \rightarrow a \{(b,c)\} \} )
\end{align*}
\noindent were we have omitted the equivalent coverings when Black starts with  $b$  or $e$. 
As a result, this configuration of 3 groups is covered with just 5 markers, a reduction of 1 node.
The diagram in Fig. \ref{fig:Triangle} shows an example position. Black has only three groups left (indicated by the dashed lines). A naive HJ-pairing fails. However, recognizing that this position contains a Triangle Configuration (for instance with $a$, $b$, $c$, $d$, and $e$ corresponding to a1, a4, d1, c1, and c2) we see that these five squares form a matching set for these three groups, proving that this position is (at most) a draw for Black.

\subsection{The Square Configuration}

\begin{figure}[!h]
\begin{center}
\begin{tikzpicture}[scale=0.7]
\draw[fill=black] (-1.5,-1.5) circle (3.0pt);
\node at (-1.8,-1.8) {$a$};
\draw[fill=black] (-1.5,1.5) circle (3.0pt);
\node at (-1.8,1.8) {$b$};
\draw[fill=black] (1.5,1.5) circle (3.0pt);
\node at (1.8,1.8) {$c$};
\draw[fill=black] (1.5,-1.5) circle (3.0pt);
\node at (1.8,-1.8) {$d$};
\draw[fill=black] (0.0,-1.5) circle (3.0pt);
\node at (0.0,-1.8) {$e$};
\draw[fill=black] (0.0,1.5) circle (3.0pt);
\node at (0.0,1.8) {$f$};
\draw[fill=black] (1.5,0.0) circle (3.0pt);
\node at (1.8,0.0) {$g$};
\draw[dashed] (-1.5,-1.5) -- (-1.5,1.5);
\draw[dashed] (-1.5,-1.5) -- (1.5,-1.5);
\draw[dashed] (-1.5,1.5) -- (1.5,1.5);
\draw[dashed] (1.5,-1.5) -- (1.5,1.5);
\end{tikzpicture} 
\hspace*{1.0cm}
\begin{tikzpicture}[scale=0.9]
\draw[fill=black] (-1.05,-1.05) circle (2.0pt);
\draw[fill=black] (-1.05,1.05) circle (2.0pt);
\draw[fill=black] (1.05,1.05) circle (2.0pt);
\draw[fill=black] (1.05,-1.05) circle (2.0pt);
\draw[fill=black] (-0.35,1.05) circle (2.0pt);
\draw[fill=black] (1.05,-0.35) circle (2.0pt);
\draw[fill=black] (-0.35,-1.05) circle (2.0pt);
\draw[dashed] (-1.05,-1.05) -- (-1.05,1.05);
\draw[dashed] (-1.05,1.05) -- (1.05,1.05);
\draw[dashed] (-1.05,-1.05) -- (1.05,-1.05);
\draw[dashed] (1.05,-1.05) -- (1.05,1.05);
\input{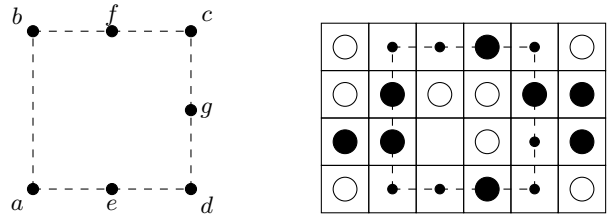}
\matrix[square matrix=.63cm]
{
	&	&	& 	& 	& 	\\
	&	&	& 	& 	& 	\\
	&	&	& 	& 	& 	\\
	&	&	& 	& 	& 	\\
};
\draw (-1.75,-1.05) circle (5.0pt);
\draw[fill=black] (0.35,-1.05) circle (5.0pt);
\draw (1.75,-1.05) circle (5.0pt);
\draw[fill=black] (-1.75,-0.35) circle (5.0pt);
\draw[fill=black] (-1.05,-0.35) circle (5.0pt);
\draw (0.35,-0.35) circle (5.0pt);
\draw[fill=black] (1.75,-0.35) circle (5.0pt);
\draw (-1.75,0.35) circle (5.0pt);
\draw[fill=black] (-1.05,0.35) circle (5.0pt);
\draw (-0.35,0.35) circle (5.0pt);
\draw (0.35,0.35) circle (5.0pt);
\draw[fill=black] (1.05,0.35) circle (5.0pt);
\draw[fill=black] (1.75,0.35) circle (5.0pt);
\draw (-1.75,1.05) circle (5.0pt);
\draw[fill=black] (0.35,1.05) circle (5.0pt);
\draw (1.75,1.05) circle (5.0pt);
\end{tikzpicture}
\end{center}
\vspace*{-0.2cm}
\caption{Sketch  and example position  illustrating the Square Configuration.}
    \label{fig:Square}
\end{figure}

We can easily define a similar configuration using a Square, with one side having two marker nodes (the {\em main markers}), the other three sides with three markers, see   Fig. \ref{fig:Square}.

\vspace*{-0.5cm}

\noindent This covering is rather similar as for the Triangle, just continuing after the first move and response the naive HJ-pairing for one more group along the cycle. Therefore, the formal notation for a matching set of this configuration is:

\vspace*{-0.3cm}
\begin{align*}
M = (\{&a \mhyphen g\}, \{ab,aed,bfc,cgd \},\\
\{ &a \rightarrow b \{(c,g),(d,e)\},  
     c \rightarrow a \{(b,f),(d,g)\}, \\ 
    &e \rightarrow a \{(b,c),(d,g)\}, 
     g \rightarrow a \{(b,f),(c,d)\}\})
\end{align*}

\noindent were we have omitted the equivalent coverings when Black starts with $b$, $d$, or $f$.
As a result, this configuration of 4 groups is covered with just 7 markers, again a reduction of 1 node.
An example position  is given in  the diagram in Fig. \ref{fig:Square}. 

\subsection{Arbitrary Cycle Configurations}

Closer examination shows that we can generalize this to any cycle of size $n$ ($C_n$), with all corners as markers and all sides but one having one additional marker. For $C_3$ and $C_4$ this yields the Triangle  and the Square Configurations above.
Note that the relative gain decreases with increasing $n$, as the number of nodes reduction per cycle is just 1, irrespective of cycle size. Of course it is possible to have multiple independent cycle configurations simultaneously, where independent means that the configurations have disjoint node sets as well as disjoint group sets (though groups of different configurations may intersect).

\section{Set Matching using BiCycle Configurations}

The cycle configurations just give a reduction of 1  node  per independent cycle. It therefore would be advantageous when multiple cycles interact in such a way that they can be used in a  more efficient matching set. For this we need overlapping node sets and possibly overlapping group sets as well. We need at least two common nodes in the node sets, since otherwise the first player  can start with the common node, after which we have two disjoint results without any reduction at all. 

In this section we will investigate configurations consisting of two interacting triangles, starting with one triangle and adding a second one. When the two triangles have  three common sides, they necessarily coincide, and no gain can be obtained.
We therefore consider the following 3 cases: 1)  two triangles with two common sides and one new side; 2) two triangles with one common side and two new sides; and 3) two triangles with no common side and three new sides.

\subsection{Configurations of Two Triangles with One New Side}

\begin{figure}[!h]
\begin{center}
\begin{tikzpicture}[scale=0.7]
\draw[fill=black] (-1.0,-1.732) circle (3.0pt);
\node at (-1.3,-1.732) {$a$};
\draw[fill=black] (-1.0,1.732) circle (3.0pt);
\node at (-1.3,1.732) {$b$};
\draw[fill=black] (2.0,0.0) circle (3.0pt);
\node at (2.3,0.0) {$c$};
\draw[fill=black] (0.5,-0.866) circle (3.0pt);
\node at (0.8,-1.266) {$d$};
\draw[fill=black] (0.5,0.866) circle (3.0pt);
\node at (0.8,1.266) {$e$};
\draw[fill=black] (0.5,0.0) circle (3.0pt);
\node at (0.8,0.0) {$f$};
\draw[dashed,ultra thick] (2.0,0.0) -- (-1.0,1.732);
\draw[dashed,ultra thick] (2.0,0.0) -- (-1.0,-1.732);
\draw[dashed] (-1.0,1.732) -- (-1.0,-1.732);
\draw[dashed] (0.5,1.732) -- (0.5,-1.732);
\end{tikzpicture} 
\hspace*{1.0cm}
\begin{tikzpicture}[scale=0.9]
\draw[fill=black] (-1.4,-1.05) circle (2.0pt);
\draw[fill=black] (-1.4,1.05) circle (2.0pt);
\draw[fill=black] (0.0,-1.05) circle (2.0pt);
\draw[fill=black] (0.0,-0.35) circle (2.0pt);
\draw[fill=black] (0.7,-1.05) circle (2.0pt);
\draw[fill=black] (0.0,1.05) circle (2.0pt);
\draw[dashed] (-1.4,-1.05) -- (-1.4,1.05);
\draw[dashed] (-1.4,1.05) -- (0.7,-1.05);
\draw[dashed] (-1.4,-1.05) -- (0.7,-1.05);
\draw[dashed] (0.0,-1.05) -- (0.0,1.05);
\input{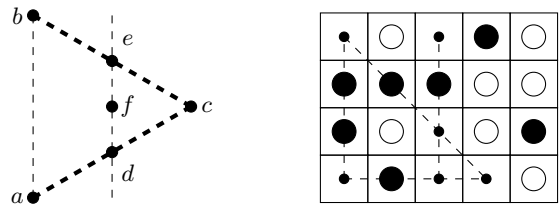}
\matrix[square matrix=.63cm]
{
	&	& 	& 	& 	\\
	& 	&	& 	& 	\\
	& 	& 	&	& 	\\
	& 	& 	& 	&	\\
};
\draw[fill=black] (-0.7,-1.05) circle (5pt);
\draw (1.4,-1.05) circle (5pt);
\draw[fill=black] (-1.4,-0.35) circle (5pt);
\draw (-0.7,-0.35) circle (5pt);
\draw[fill=black] (1.4,-0.35) circle (5pt);
\draw[fill=black] (-1.4,0.35) circle (5pt);
\draw[fill=black] (-0.7,0.35) circle (5pt);
\draw[fill=black] (0.0,0.35) circle (5pt);
\draw (0.7,0.35) circle (5pt);
\draw (1.4,0.35) circle (5pt);
\draw (-0.7,1.05) circle (5pt);
\draw (0.7,-0.35) circle (5pt);
\draw[fill=black] (0.7,1.05) circle (5pt);
\draw (1.4,1.05) circle (5pt);
\end{tikzpicture}
\end{center}
\vspace*{-0.2cm}
\caption{Sketch   and example position   illustrating a Triangle/Line Configuration.}
    \label{fig:TriangleLine}
\end{figure}

When  the two triangles have two common sides and one new side we have several possibilities for the new side. We just give one case, namely where the additional side uses the two edge markers of the first triangle. Since this configuration can be equally well seen as a cycle with an additional line with additional marker intersecting the two edge markers of the triangle, we denote this configuration as a Triangle/Line Configuration, see Fig. \ref{fig:TriangleLine}. The thick dashed lines denote the common sides of the triangles.

\vspace*{-0.3cm}

In this configuration, the main markers are $a$ and $b$, like in the simple Triangle Configuration with triangle $abc$. Since depending on the moves played in this triangle either $d$ or $e$ will not be used for covering the sides of the first triangle, they may be used profitably with one additional edge marker for a new side {\em dfe} of the second triangle $dec$.
The formal notation for a matching set of this configuration is:

\vspace*{-0.3cm}
\begin{align*}
M = (\{&a \mhyphen f\}, \{ab,adc,bec,dfe\}, \\
\{&a \rightarrow b\{(c,d),(e,f)\}, 
    c \rightarrow a\{(b,e),(d,f)\}, \\
   &d \rightarrow a\{(b,c),(e,f)\}, 
    f \rightarrow a\{(b,c),(d,e)\}\})
\end{align*}

\noindent were we have omitted the equivalent coverings when Black starts with $b$ or $e$.
As a result, this configuration of 4 groups is covered with just 6 markers, a reduction of 2 nodes.
An example position  is given in  the diagram in Fig. \ref{fig:TriangleLine}. Note that edge markers need not be physically situated between the corner markers in a group.

\begin{figure}[!h]
\begin{center}
\begin{tikzpicture}[scale=0.7]
\draw[fill=black] (-1.5,-1.5) circle (3.0pt);
\node at (-1.8,-1.8) {$a$};
\draw[fill=black] (-1.5,1.5) circle (3.0pt);
\node at (-1.8,1.8) {$b$};
\draw[fill=black] (1.5,1.5) circle (3.0pt);
\node at (1.8,1.8) {$c$};
\draw[fill=black] (1.5,-1.5) circle (3.0pt);
\node at (1.8,-1.8) {$d$};
\draw[fill=black] (0.0,-1.5) circle (3.0pt);
\node at (0.0,-1.8) {$e$};
\draw[fill=black] (0.0,1.5) circle (3.0pt);
\node at (0.0,1.8) {$f$};
\draw[fill=black] (1.5,0.0) circle (3.0pt);
\node at (1.8,0.0) {$g$};
\draw[fill=black] (0.0,0.0) circle (2.0pt);
\node at (0.3,0.0) {$h$};
\draw[dashed] (-1.5,-1.5) -- (-1.5,1.5);
\draw[dashed,ultra thick] (-1.5,-1.5) -- (1.5,-1.5);
\draw[dashed,ultra thick] (-1.5,1.5) -- (1.5,1.5);
\draw[dashed,ultra thick] (1.5,-1.5) -- (1.5,1.5);
\draw[dashed] (0.0,-1.5) -- (0.0,1.5);
\end{tikzpicture} 
\hspace*{1.0cm}
\begin{tikzpicture}[scale=0.9]
\draw[fill=black] (-1.05,-1.05) circle (2.0pt);
\draw[fill=black] (-1.05,1.05) circle (2.0pt);
\draw[fill=black] (1.05,1.05) circle (2.0pt);
\draw[fill=black] (1.05,-1.05) circle (2.0pt);
\draw[fill=black] (-0.35,1.05) circle (2.0pt);
\draw[fill=black] (1.05,-0.35) circle (2.0pt);
\draw[fill=black] (-0.35,-1.05) circle (2.0pt);
\draw[fill=black] (-0.35,-0.35) circle (2.0pt);
\draw[dashed] (-1.05,-1.05) -- (-1.05,1.05);
\draw[dashed] (-1.05,1.05) -- (1.05,1.05);
\draw[dashed] (-1.05,-1.05) -- (1.05,-1.05);
\draw[dashed] (1.05,-1.05) -- (1.05,1.05);
\draw[dashed] (-0.35,-1.05) -- (-0.35,1.05);
\input{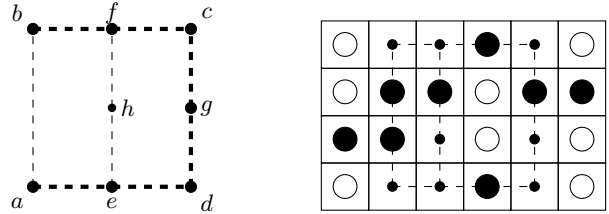}
\matrix[square matrix=.63cm]
{
	&	&	& 	& 	& 	\\
	&	&	& 	& 	& 	\\
	&	&	& 	& 	& 	\\
	&	&	& 	& 	& 	\\
};
\draw (-1.75,-1.05) circle (5.0pt);
\draw[fill=black] (0.35,-1.05) circle (5.0pt);
\draw (1.75,-1.05) circle (5.0pt);
\draw[fill=black] (-1.75,-0.35) circle (5.0pt);
\draw[fill=black] (-1.05,-0.35) circle (5.0pt);
\draw (0.35,-0.35) circle (5.0pt);
\draw (1.75,-0.35) circle (5.0pt);
\draw (-1.75,0.35) circle (5.0pt);
\draw[fill=black] (-1.05,0.35) circle (5.0pt);
\draw[fill=black] (-0.35,0.35) circle (5.0pt);
\draw (0.35,0.35) circle (5.0pt);
\draw[fill=black] (1.05,0.35) circle (5.0pt);
\draw[fill=black] (1.75,0.35) circle (5.0pt);
\draw (-1.75,1.05) circle (5.0pt);
\draw[fill=black] (0.35,1.05) circle (5.0pt);
\draw (1.75,1.05) circle (5.0pt);
\end{tikzpicture}
\end{center}
\vspace*{-0.2cm}
\caption{Sketch and example position   illustrating the Square/Line Configuration.}
    \label{fig:SquareLine}
\end{figure}

The idea of two cycles with one new side also works with other cycles. We just give one example, see Fig. \ref{fig:SquareLine}.

\vspace*{-0.3cm}

\noindent
Here we have two overlapping Squares with three common sides (a Square/Line Configuration), where the second Square has one new side  with one additional marker. The formal notation for a matching set of this configuration is:

\vspace*{-0.3cm}
\begin{align*}
M = (\{&a \mhyphen h\}, \{ab,aed,bfc,cgd,ehf \},\\
\{ &a \rightarrow b \{(c,g),(d,e), (f,h)\},  
     c \rightarrow a \{(b,f),(d,g),(e,h)\},  \\
    &e \rightarrow a \{(b,c),(d,g),(f,h)\}, 
      g \rightarrow a \{(b,f),(c,d),(e,h)\}, \\
    &h \rightarrow a \{(b,c),(d,g),(e,f)\}\})
\end{align*}

\noindent were we have omitted the equivalent coverings when Black starts with $b$, $d$, or $f$.
As a result, this configuration of 5 groups is covered with just 8 markers, again a reduction of 2 nodes. The diagram in Fig. \ref{fig:SquareLine} shows an example position  using this configuration.

\subsection{Configurations of Two Triangles with Two New Sides}

\begin{figure}[!h]
\begin{center}
\begin{tikzpicture}[scale=0.7]
\draw[fill=black] (-1.0,-1.732) circle (2.0pt);
\node at (-1.0,-2.032) {$a$};
\draw[fill=black] (-1.0,1.732) circle (2.0pt);
\node at (-1.0,2.032) {$b$};
\draw[fill=black] (2.0,0.0) circle (2.0pt);
\node at (2.3,0.0) {$c$};
\draw[fill=black] (0.5,-0.866) circle (2.0pt);
\node at (0.65,-1.166) {$d$};
\draw[fill=black] (0.5,0.866) circle (2.0pt);
\node at (0.65,1.166) {$e$};
\draw[fill=black] (-4.0,0.0) circle (2.0pt);
\node at (-4.3,0.0) {$f$};
\draw[fill=black] (-2.5,-0.866) circle (2.0pt);
\node at (-2.65,-1.166) {$g$};
\draw[fill=black] (-2.5,0.866) circle (2.0pt);
\node at (-2.65,1.166) {$h$};
\draw[dashed] (2.0,0.0) -- (-1.0,1.732);
\draw[dashed] (2.0,0.0) -- (-1.0,-1.732);
\draw[dashed, ultra thick] (-1.0,1.732) -- (-1.0,-1.732);
\draw[dashed] (-4.0,0.0) -- (-1.0,1.732);
\draw[dashed] (-4.0,0.0) -- (-1.0,-1.732);
\end{tikzpicture} 
\hspace*{1.0cm}
\begin{tikzpicture}[scale=0.9]
\draw[fill=black] (0.0,-1.05) circle (2.0pt);
\draw[fill=black] (0.0,1.05) circle (2.0pt);
\draw[fill=black] (1.4,0.35) circle (2.0pt);
\draw[fill=black] (1.4,1.05) circle (2.0pt);
\draw[fill=black] (2.1,1.05) circle (2.0pt);
\draw[fill=black] (-1.4,-0.35) circle (2.0pt);
\draw[fill=black] (-1.4,-1.05) circle (2.0pt);
\draw[fill=black] (-2.1,-1.05) circle (2.0pt);
\draw[dashed] (0.0,-1.05) -- (0.0,1.05);
\draw[dashed] (0.0,1.05) -- (2.1,1.05);
\draw[dashed] (0.0,-1.05) -- (2.1,1.05);
\draw[dashed] (0.0,1.05) -- (-2.1,-1.05);
\draw[dashed] (0.0,-1.05) -- (-2.1,-1.05);
\input{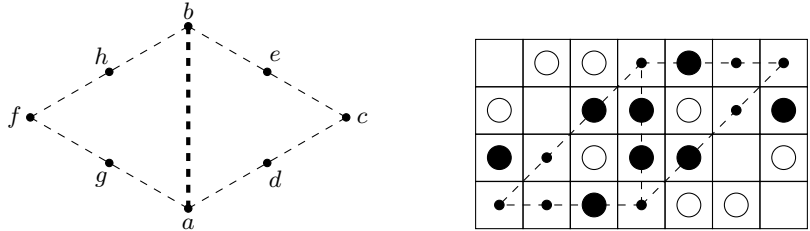}
\matrix[square matrix=.63cm]
{
	&	&	&	& 	& 	& 	\\
	&	&	&	& 	& 	& 	\\
	&	&	&	& 	& 	& 	\\
	&	&	&	& 	& 	& 	\\
};
\draw[fill=black] (-0.7,-1.05) circle (5.0pt);
\draw (0.7,-1.05) circle (5.0pt);
\draw (1.4,-1.05) circle (5.0pt);
\draw[fill=black] (-2.1,-0.35) circle (5.0pt);
\draw (-0.7,-0.35) circle (5.0pt);
\draw[fill=black] (0.0,-0.35) circle (5.0pt);
\draw[fill=black] (0.7,-0.35) circle (5.0pt);
\draw (2.1,-0.35) circle (5.0pt);
\draw[fill=black] (-0.7,0.35) circle (5.0pt);
\draw (-2.1,0.35) circle (5.0pt);
\draw[fill=black] (0.0,0.35) circle (5.0pt);
\draw (0.7,0.35) circle (5.0pt);
\draw[fill=black] (2.1,0.35) circle (5.0pt);
\draw (-1.4,1.05) circle (5.0pt);
\draw (-0.7,1.05) circle (5.0pt);
\draw[fill=black] (0.7,1.05) circle (5.0pt);
\end{tikzpicture}
\end{center}
\vspace*{-0.2cm}
\caption{Sketch   and example position   illustrating the BiTriangle Configuration.}
    \label{fig:BiTriangle}
\end{figure}

When the two triangles have one side in common, the common side should be the  side with just two markers  of the first triangle. If not, the first player could choose one of the main markers such that the second player has to respond outside the second triangle, which does not give a profitable situation. When the two triangles share this side, we have two subcases. The first is when no other nodes are shared, leading to a BiTriangle Configuration, see Fig. \ref{fig:BiTriangle}. Responses by White in these two cycles are related by mirror  symmetry.

\vspace*{-0.3cm}

\vspace*{-0.3cm}
\noindent
The formal notation for a matching set of this configuration is:
\begin{align*}
M = (&\{a \mhyphen h\}, \{ ab,adc,agf,bec,bhf\}, \\
&\{ a \rightarrow b\{(c,d),(f,g)\}, c \rightarrow a\{(b,e),(f,h)\}, d \rightarrow a\{(b,c),(f,h)\}\}) 
\end{align*}

\noindent were we have omitted the equivalent coverings when Black starts with $b$, $e$, $f$, $g$, or $h$. 
In this configuration 8 markers cover 5 groups, a reduction of 2 nodes. 
An example  position is given in the diagram in Fig. \ref{fig:BiTriangle}. 

\begin{figure}[!h]
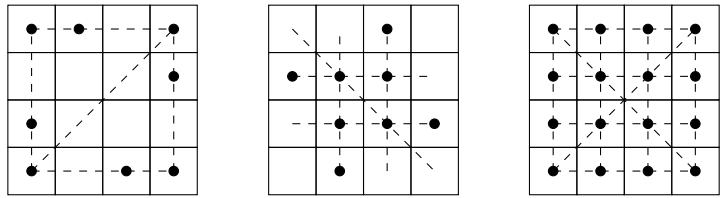

\begin{center}
\begin{tikzpicture}[scale=0.9]
\draw[fill=black] (-1.05,-1.05) circle (2.0pt);
\draw[fill=black] (-1.05,1.05) circle (2.0pt);
\draw[fill=black] (1.05,1.05) circle (2.0pt);
\draw[fill=black] (-1.05,-0.35) circle (2.0pt);
\draw[fill=black] (-0.35,1.05) circle (2.0pt);
\draw[fill=black] (0.35,-1.05) circle (2.0pt);
\draw[fill=black] (1.05,-1.05) circle (2.0pt);
\draw[fill=black] (1.05,0.35) circle (2.0pt);
\draw[dashed] (-1.05,-1.05) -- (-1.05,1.05);
\draw[dashed] (-1.05,1.05) -- (1.05,1.05);
\draw[dashed] (-1.05,-1.05) -- (1.05,1.05);
\draw[dashed] (-1.05,-1.05) -- (1.05,-1.05);
\draw[dashed] (1.05,-1.05) -- (1.05,1.05);
\input{boardintro.tex}
\matrix[square matrix=.63cm]
{
	& 	& 	& 	\\
	& 	& 	& 	\\
	& 	& 	& 	\\
	& 	& 	& 	\\
};
\end{tikzpicture}
\hspace*{0.5cm}
\begin{tikzpicture}[scale=0.9]
\draw[fill=black] (-0.35,-0.35) circle (2.0pt);
\draw[fill=black] (-0.35,0.35) circle (2.0pt);
\draw[fill=black] (0.35,0.35) circle (2.0pt);
\draw[fill=black] (0.35,-0.35) circle (2.0pt);
\draw[fill=black] (-0.35,-1.05) circle (2.0pt);
\draw[fill=black] (0.35,1.05) circle (2.0pt);
\draw[fill=black] (-1.05,0.35) circle (2.0pt);
\draw[fill=black] (1.05,-0.35) circle (2.0pt);
\draw[dashed] (-0.35,-1.05) -- (-0.35,1.05);
\draw[dashed] (-1.05,0.35) -- (1.05,0.35);
\draw[dashed] (-1.05,1.05) -- (1.05,-1.05);
\draw[dashed] (-1.05,-0.35) -- (1.05,-0.35);
\draw[dashed] (0.35,-1.05) -- (0.35,1.05);
\input{boardintro.tex}
\matrix[square matrix=.63cm]
{
	& 	& 	& 	\\
	& 	& 	& 	\\
	& 	& 	& 	\\
	& 	& 	& 	\\
};
\end{tikzpicture}
\hspace*{0.5cm}
\begin{tikzpicture}[scale=0.9]
\draw[fill=black] (-1.05,-1.05) circle (2.0pt);
\draw[fill=black] (-1.05,1.05) circle (2.0pt);
\draw[fill=black] (1.05,1.05) circle (2.0pt);
\draw[fill=black] (-1.05,-0.35) circle (2.0pt);
\draw[fill=black] (-0.35,1.05) circle (2.0pt);
\draw[fill=black] (-0.35,-1.05) circle (2.0pt);
\draw[fill=black] (1.05,-1.05) circle (2.0pt);
\draw[fill=black] (1.05,-0.35) circle (2.0pt);
\draw[fill=black] (-0.35,-0.35) circle (2.0pt);
\draw[fill=black] (-0.35,0.35) circle (2.0pt);
\draw[fill=black] (0.35,0.35) circle (2.0pt);
\draw[fill=black] (0.35,-0.35) circle (2.0pt);
\draw[fill=black] (0.35,-1.05) circle (2.0pt);
\draw[fill=black] (0.35,1.05) circle (2.0pt);
\draw[fill=black] (-1.05,0.35) circle (2.0pt);
\draw[fill=black] (1.05,0.35) circle (2.0pt);
\draw[dashed] (-1.05,-1.05) -- (-1.05,1.05);
\draw[dashed] (-1.05,1.05) -- (1.05,1.05);
\draw[dashed] (-1.05,-1.05) -- (1.05,1.05);
\draw[dashed] (-1.05,-1.05) -- (1.05,-1.05);
\draw[dashed] (1.05,-1.05) -- (1.05,1.05);
\draw[dashed] (-0.35,-1.05) -- (-0.35,1.05);
\draw[dashed] (-1.05,0.35) -- (1.05,0.35);
\draw[dashed] (-1.05,1.05) -- (1.05,-1.05);
\draw[dashed] (-1.05,-0.35) -- (1.05,-0.35);
\draw[dashed] (0.35,-1.05) -- (0.35,1.05);
\input{boardintro.tex}
\matrix[square matrix=.63cm]
{
	& 	& 	& 	\\
	& 	& 	& 	\\
	& 	& 	& 	\\
	& 	& 	& 	\\
};
\end{tikzpicture}
\end{center}
\vspace*{-0.2cm}
\caption{Example positions illustrating two BiTriangle Configurations (left and middle) and their combination (right) on the empty $4 \times 4$ board.}
    \label{fig:Empty4x4}
\end{figure}

Note that  several independent configurations may be used simultaneously. As a notorious example, considering the empty $4 \times 4$ board, we  can apply a  BiTriangle Configuration as in Fig. \ref{fig:Empty4x4} (left) or alternatively apply another  BiTriangle Configuration as in Fig. \ref{fig:Empty4x4} (middle).
 
\vspace*{-0.3cm}

Since these two matching sets do not interact, we may combine the two configurations for the empty $4 \times 4$ board (they  just fit ``into'' each other), resulting in Fig. \ref{fig:Empty4x4} (right). As a consequence all 10  groups  are matched by the 16 empty squares. So this new pairing method proves that the $4 \times 4$ board is a draw, investigating just 1 node (a kind of perfect solving \cite{Uiterwijk2014}).

\begin{figure}[!h]
\begin{center}
\begin{tikzpicture}[scale=0.7]
\draw[fill=black] (-1.0,-1.732) circle (2.0pt);
\node at (-1.0,-2.032) {$a$};
\draw[fill=black] (-1.0,1.732) circle (2.0pt);
\node at (-1.0,2.032) {$b$};
\draw[fill=black] (2.0,-1.732) circle (2.0pt);
\node at (2.3,-2.032) {$c$};
\draw[fill=black] (0.5,-1.732) circle (2.0pt);
\node at (0.5,-2.032) {$d$};
\draw[fill=black] (0.5,0.0) circle (2.0pt);
\node at (0.8,0.0) {$e$};
\draw[fill=black] (2.0,1.732) circle (2.0pt);
\node at (2.3,2.032) {$f$};
\draw[fill=black] (0.5,1.732) circle (2.0pt);
\node at (0.5,2.032) {$g$};
\draw[dashed] (2.0,-1.732) -- (-1.0,-1.732); 
\draw[dashed] (2.0,-1.732) -- (-1.0,1.732); 
\draw[dashed, ultra thick] (-1.0,1.732) -- (-1.0,-1.732); 
\draw[dashed] (2.0,1.732) -- (-1.0,-1.732); 
\draw[dashed] (2.0,1.732) -- (-1.0,1.732); 
\end{tikzpicture} 
\hspace*{1.0cm}
\begin{tikzpicture}[scale=0.9]
\draw[fill=black] (-0.7,-1.05) circle (2.0pt); 
\draw[fill=black] (-0.7,0.35) circle (2.0pt); 
\draw[fill=black] (0.7,0.35) circle (2.0pt); 
\draw[fill=black] (0.0,-0.35) circle (2.0pt); 
\draw[fill=black] (0.0,0.35) circle (2.0pt); 
\draw[fill=black] (0.0,-1.05) circle (2.0pt); 
\draw[fill=black] (0.7,-1.05) circle (2.0pt); 
\draw[dashed] (-0.7,-1.05) -- (-0.7,1.05);
\draw[dashed] (-0.7,0.35) -- (1.4,0.35);
\draw[dashed] (-0.7,-1.05) -- (1.4,1.05);
\draw[dashed] (-0.7,-1.05) -- (1.4,-1.05);
\draw[dashed] (0.7,-1.05) -- (-1.4,1.05);
\input{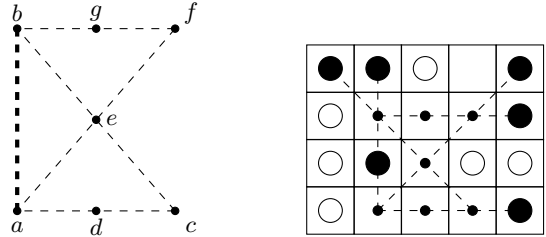}
\matrix[square matrix=.63cm]
{
	&	& 	& 	& 	\\
	& 	&	& 	& 	\\
	& 	& 	&	& 	\\
	& 	& 	& 	&	\\
};
\draw  (-1.4,-1.05) circle (5.0pt); 
\draw[fill=black] (1.4,-1.05) circle (5.0pt); 
\draw  (-1.4,-0.35) circle (5.0pt); 
\draw[fill=black] (-0.7,-0.35) circle (5.0pt); 
\draw  (0.7,-0.35) circle (5.0pt); 
\draw  (1.4,-0.35) circle (5.0pt); 
\draw  (-1.4,0.35) circle (5.0pt); 
\draw[fill=black]  (1.4,0.35) circle (5.0pt); 
\draw[fill=black]  (-1.4,1.05) circle (5.0pt); 
\draw[fill=black]  (-0.7,1.05) circle (5.0pt); 
\draw  (0.0,1.05) circle (5.0pt); 
\draw[fill=black]  (1.4,1.05) circle (5.0pt); 
\end{tikzpicture}
\end{center}
\vspace*{-0.2cm}
\caption{Sketch  and example position   illustrating the BiTriangleX Configuration.}
    \label{fig:BiTriangleX}
\end{figure}

As the second subcase we consider the possibility  that the two triangles  also share another edge marker, see Fig. \ref{fig:BiTriangleX}. We denote this configuration as the BiTriangleX Configuration.

\vspace*{-0.3cm}

\noindent
The formal notation for a matching set of this configuration  is:

\vspace*{-0.3cm}
\begin{align*}
M = (\{&a \mhyphen g\}, \{ ab,adc,aef,bec,bgf\}, \\
\{ &a \rightarrow b\{(c,d),(e,f)\}, 
      c \rightarrow a\{(b,e),(f,g)\}, \\
   &d \rightarrow a\{(b,f),(c,e)\}, 
     e \rightarrow a\{(b,c),(f,g)\}\}) 
\end{align*}

\noindent were we have omitted the equivalent coverings when Black starts with  $b$, $f$, or $g$. 
In this configuration 7 markers  cover 5 groups, a reduction of 3 nodes. 
An example position is given in the diagram  in Fig. \ref{fig:BiTriangleX}. 


\subsection{Configurations of Two Triangles with Three New Sides}

\begin{figure}[!h]
\begin{center}
\begin{tikzpicture}[scale=0.7]
\draw[fill=black] (-1.0,-1.732) circle (2.0pt);
\node at (-1.3,-1.732) {$a$};
\draw[fill=black] (-1.0,0.0) circle (2.0pt);
\node at (-1.3,0.0) {$b$};
\draw[fill=black] (-1.0,1.732) circle (2.0pt);
\node at (-1.3,1.732) {$c$};
\draw[fill=black] (0.5,-0.866) circle (2.0pt);
\node at (0.5,-1.166) {$d$};
\draw[fill=black] (0.5,0.866) circle (2.0pt);
\node at (0.5,1.166) {$e$};
\draw[fill=black] (2.0,-1.732) circle (2.0pt);
\node at (2.3,-1.732) {$f$};
\draw[fill=black] (2.0,0.0) circle (2.0pt);
\node at (2.3,0.0) {$g$};
\draw[fill=black] (2.0,1.732) circle (2.0pt);
\node at (2.3,1.732) {$h$};
\draw[dashed] (2.0,0.0) -- (-1.0,1.732);
\draw[dashed] (2.0,0.0) -- (-1.0,-1.732);
\draw[dashed] (-1.0,1.732) -- (-1.0,-1.732);
\draw[dashed] (-1.0,0.0) -- (2.0,-1.732);
\draw[dashed] (-1.0,0.0) -- (2.0,1.732);
\draw[dashed] (2.0,-1.732) -- (2.0,1.732);
\end{tikzpicture} 
\hspace*{1.0cm}
\begin{tikzpicture}[scale=0.9]
\draw[fill=black] (-1.05,-0.35) circle (2.0pt);
\draw[fill=black] (-0.35,-1.05) circle (2.0pt);
\draw[fill=black] (-0.35,-0.35) circle (2.0pt);
\draw[fill=black] (-0.35,0.35) circle (2.0pt);
\draw[fill=black] (0.35,-0.35) circle (2.0pt);
\draw[fill=black] (0.35,0.35) circle (2.0pt);
\draw[fill=black] (0.35,1.05) circle (2.0pt);
\draw[fill=black] (1.05,0.35) circle (2.0pt);
\draw[dashed] (-1.05,-0.35) -- (1.05,-0.35);
\draw[dashed] (-1.75,-1.05) -- (0.35,1.05);
\draw[dashed] (0.35,-1.05) -- (0.35,1.05);
\draw[dashed] (-0.35,-1.05) -- (-0.35,1.05);
\draw[dashed] (-1.05,0.35) -- (1.05,0.35);
\draw[dashed] (-0.35,-1.05) -- (1.75,1.05);
\input{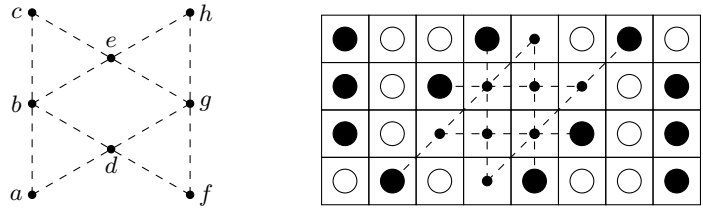}
\matrix[square matrix=.63cm]
{
	&	&	&	&	& 	& 	& 	\\
	&	&	&	&	& 	& 	& 	\\
	&	&	&	&	& 	& 	& 	\\
	&	&	&	&	& 	& 	& 	\\
};
\draw[fill=black] (-1.75,-1.05) circle (5.0pt);
\draw (-1.05,-1.05) circle (5.0pt);
\draw[fill=black] (0.35,-1.05) circle (5.0pt);
\draw (1.05,-1.05) circle (5.0pt);
\draw (1.75,-1.05) circle (5.0pt);
\draw (-1.75,-0.35) circle (5.0pt);
\draw[fill=black] (1.05,-0.35) circle (5.0pt);
\draw (1.75,-0.35) circle (5.0pt);
\draw (-1.75,0.35) circle (5.0pt);
\draw[fill=black] (-1.05,0.35) circle (5.0pt);
\draw (1.75,0.35) circle (5.0pt);
\draw (-1.75,1.05) circle (5.0pt);
\draw (-1.05,1.05) circle (5.0pt);
\draw[fill=black] (-0.35,1.05) circle (5.0pt);
\draw (1.05,1.05) circle (5.0pt);
\draw[fill=black] (1.75,1.05) circle (5.0pt);
\draw[fill=black] (-2.45,-0.35) circle (5.0pt);
\draw[fill=black] (-2.45,0.35) circle (5.0pt);
\draw[fill=black] (2.45,-0.35) circle (5.0pt);
\draw[fill=black] (2.45,0.35) circle (5.0pt);
\draw[fill=black] (-2.45,1.05) circle (5.0pt);
\draw (-2.45,-1.05) circle (5.0pt);
\draw[fill=black] (2.45,-1.05) circle (5.0pt);
\draw (2.45,1.05) circle (5.0pt);
\end{tikzpicture}
\end{center}
\vspace*{-0.2cm}
\caption{Sketch  and example position  illustrating the FlatStar Configuration.}
    \label{fig:FlatStar}
\end{figure}

Another idea is to use the triangles without reduction, but sharing common markers. There are many ways to do this. We just give one profitable case, see Fig. \ref{fig:FlatStar}. We call this configuration the FlatStar Configuration.

\vspace*{-0.3cm}

In this configuration the pairing set consists of the eight nodes $a$--$h$. 
It is easy to verify that after every first marked node played by the first player, the second player has a suitable response by $b$ or $g$ (the two main markers), leaving a situation were all remaining  groups can be covered with a naive HJ-pairing.
As a result  the pairing set  $a$-$h$ covers all 6 groups.
The formal notation for a matching set of this configuration is:

\vspace*{-0.3cm}
\begin{align*}
M = (\{&a \mhyphen h\}, \{abc,adg,bdf,beh,ceg,fgh\}, \\ 
\{&a \rightarrow b \{(c,e)(d,g)(f,h)\}, 
     b \rightarrow g \{(a,c)(d,f)(e,h)\}, \\
   &d \rightarrow b \{(a,g)(c,e)(f,h)\}\}) 
\end{align*}

\noindent were we have omitted the equivalent coverings when Black starts with $c$, $e$, $f$, $g$, or $h$.
This configuration is particularly efficient, needing just 8 markers for 6 groups, a reduction of 4 nodes.
An example position is given in the diagram  in Fig. \ref{fig:FlatStar}. 

\section{Set Matching using PolyCycle Configurations}

When more than two cycles are involved (further denoted as PolyCycles), the number of possible cases grows considerably. We just give two possibilities. The first possibility is to form an additional cycle  using all-but-one common sides with one cycle, leading to PolyCycle/Line Configurations, see Section \ref{subsec:PolyCycleLine}. The second possibility is when more than two cycles share one common edge, see Section \ref{subsec:PolyCommon}.

\subsection{PolyCycle/Line Configurations \label{subsec:PolyCycleLine}}

\begin{figure}[!h]
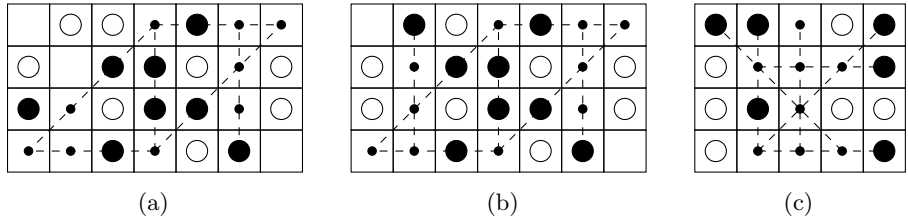

\begin{center}
\begin{tikzpicture}[scale=0.8]
\draw[fill=black] (0.0,-1.05) circle (2.0pt);
\draw[fill=black] (0.0,1.05) circle (2.0pt);
\draw[fill=black] (1.4,0.35) circle (2.0pt);
\draw[fill=black] (1.4,1.05) circle (2.0pt);
\draw[fill=black] (2.1,1.05) circle (2.0pt);
\draw[fill=black] (-1.4,-0.35) circle (2.0pt);
\draw[fill=black] (-1.4,-1.05) circle (2.0pt);
\draw[fill=black] (-2.1,-1.05) circle (2.0pt);
\draw[fill=black] (1.4,-0.35) circle (2.0pt);
\draw[dashed] (0.0,-1.05) -- (0.0,1.05);
\draw[dashed] (0.0,1.05) -- (2.1,1.05);
\draw[dashed] (0.0,-1.05) -- (2.1,1.05);
\draw[dashed] (0.0,1.05) -- (-2.1,-1.05);
\draw[dashed] (0.0,-1.05) -- (-2.1,-1.05);
\draw[dashed] (1.4,-1.05) -- (1.4,1.05);
\input{boardintro.tex}
\matrix[square matrix=.56cm]
{
	&	&	&	& 	& 	& 	\\
	&	&	&	& 	& 	& 	\\
	&	&	&	& 	& 	& 	\\
	&	&	&	& 	& 	& 	\\
};
\draw[fill=black] (-0.7,-1.05) circle (5.0pt);
\draw (0.7,-1.05) circle (5.0pt);
\draw[fill=black] (1.4,-1.05) circle (5.0pt);
\draw[fill=black] (-2.1,-0.35) circle (5.0pt);
\draw (-0.7,-0.35) circle (5.0pt);
\draw[fill=black] (0.0,-0.35) circle (5.0pt);
\draw[fill=black] (0.7,-0.35) circle (5.0pt);
\draw (2.1,-0.35) circle (5.0pt);
\draw[fill=black] (-0.7,0.35) circle (5.0pt);
\draw (-2.1,0.35) circle (5.0pt);
\draw[fill=black] (0.0,0.35) circle (5.0pt);
\draw (0.7,0.35) circle (5.0pt);
\draw (2.1,0.35) circle (5.0pt);
\draw (-1.4,1.05) circle (5.0pt);
\draw (-0.7,1.05) circle (5.0pt);
\draw[fill=black] (0.7,1.05) circle (5.0pt);
\end{tikzpicture}
\hspace*{0.2cm}
\begin{tikzpicture}[scale=0.8]
\draw[fill=black] (0.0,-1.05) circle (2.0pt);
\draw[fill=black] (0.0,1.05) circle (2.0pt);
\draw[fill=black] (1.4,0.35) circle (2.0pt);
\draw[fill=black] (1.4,1.05) circle (2.0pt);
\draw[fill=black] (2.1,1.05) circle (2.0pt);
\draw[fill=black] (-1.4,-0.35) circle (2.0pt);
\draw[fill=black] (-1.4,-1.05) circle (2.0pt);
\draw[fill=black] (-2.1,-1.05) circle (2.0pt);
\draw[fill=black] (1.4,-0.35) circle (2.0pt);
\draw[fill=black] (-1.4,0.35) circle (2.0pt);
\draw[dashed] (0.0,-1.05) -- (0.0,1.05);
\draw[dashed] (0.0,1.05) -- (2.1,1.05);
\draw[dashed] (0.0,-1.05) -- (2.1,1.05);
\draw[dashed] (0.0,1.05) -- (-2.1,-1.05);
\draw[dashed] (0.0,-1.05) -- (-2.1,-1.05);
\draw[dashed] (1.4,-1.05) -- (1.4,1.05);
\draw[dashed] (-1.4,-1.05) -- (-1.4,1.05);
\input{boardintro.tex}
\matrix[square matrix=.56cm]
{
	&	&	&	& 	& 	& 	\\
	&	&	&	& 	& 	& 	\\
	&	&	&	& 	& 	& 	\\
	&	&	&	& 	& 	& 	\\
};
\draw[fill=black] (-0.7,-1.05) circle (5.0pt);
\draw (0.7,-1.05) circle (5.0pt);
\draw[fill=black] (1.4,-1.05) circle (5.0pt);
\draw  (-2.1,-0.35) circle (5.0pt);
\draw (-0.7,-0.35) circle (5.0pt);
\draw[fill=black] (0.0,-0.35) circle (5.0pt);
\draw[fill=black] (0.7,-0.35) circle (5.0pt);
\draw (2.1,-0.35) circle (5.0pt);
\draw[fill=black] (-0.7,0.35) circle (5.0pt);
\draw (-2.1,0.35) circle (5.0pt);
\draw[fill=black] (0.0,0.35) circle (5.0pt);
\draw (0.7,0.35) circle (5.0pt);
\draw (2.1,0.35) circle (5.0pt);
\draw[fill=black]  (-1.4,1.05) circle (5.0pt);
\draw (-0.7,1.05) circle (5.0pt);
\draw[fill=black] (0.7,1.05) circle (5.0pt);
\end{tikzpicture}
\hspace*{0.2cm}
\begin{tikzpicture}[scale=0.8]
\draw[fill=black] (-0.7,-1.05) circle (2.0pt); 
\draw[fill=black] (-0.7,0.35) circle (2.0pt); 
\draw[fill=black] (0.7,0.35) circle (2.0pt); 
\draw[fill=black] (0.0,-0.35) circle (2.0pt); 
\draw[fill=black] (0.0,0.35) circle (2.0pt); 
\draw[fill=black] (0.0,-1.05) circle (2.0pt); 
\draw[fill=black] (0.7,-1.05) circle (2.0pt); 
\draw[fill=black] (0.0,1.05) circle (2.0pt); 
\draw[dashed] (-0.7,-1.05) -- (-0.7,1.05);
\draw[dashed] (-0.7,0.35) -- (1.4,0.35);
\draw[dashed] (-0.7,-1.05) -- (1.4,1.05);
\draw[dashed] (-0.7,-1.05) -- (1.4,-1.05);
\draw[dashed] (0.7,-1.05) -- (-1.4,1.05);
\draw[dashed] (0.0,-1.05) -- (0.0,1.05);
\input{boardintro.tex}
\matrix[square matrix=.56cm]
{
	&	& 	& 	& 	\\
	& 	&	& 	& 	\\
	& 	& 	&	& 	\\
	& 	& 	& 	&	\\
};
\draw  (-1.4,-1.05) circle (5.0pt); 
\draw[fill=black] (1.4,-1.05) circle (5.0pt); 
\draw  (-1.4,-0.35) circle (5.0pt); 
\draw[fill=black] (-0.7,-0.35) circle (5.0pt); 
\draw  (0.7,-0.35) circle (5.0pt); 
\draw  (1.4,-0.35) circle (5.0pt); 
\draw  (-1.4,0.35) circle (5.0pt); 
\draw[fill=black]  (1.4,0.35) circle (5.0pt); 
\draw[fill=black]  (-1.4,1.05) circle (5.0pt); 
\draw[fill=black]  (-0.7,1.05) circle (5.0pt); 
\draw  (0.7,1.05) circle (5.0pt); 
\draw[fill=black]  (1.4,1.05) circle (5.0pt); 
\end{tikzpicture}
\hspace*{0.4cm} (a) \hspace*{4.0cm} (b) \hspace*{3.3cm} (c)
\end{center}
\vspace*{-0.2cm}
\caption{Example positions   illustrating the BiTriangle/Line (a),  BiTriangle/BiLine (b) and BiTriangleX/Line (c) Configurations.}
    \label{fig:BiTriangleLines}
\end{figure}

If the configuration permits it may be possible to combine PolyCycles with Lines. Again, there are many ways to do this. A first possibility is shown  using a BiTriangle. From Fig. \ref{fig:BiTriangle}  and the accompanying matching set we see that in all move sequences in this configuration node $d$ or $e$ and node $g$ or $h$ are not used. This gives the opportunity for adding one more group through $d$ and $e$ with one additional edge marker  and/or adding one more group through $g$ and $h$ with one additional edge marker. This leads to BiTriangle/Line and BiTriangle/BiLine Configurations, see Fig. \ref{fig:BiTriangleLines} (a) and (b). 

\vspace*{-0.3cm}

\noindent
The first configuration needs 9 markers to cover 6 groups,  a reduction of 3 nodes; the second needs 10 markers to cover 7 groups,  a reduction of 4 nodes. Since these configurations are evident from the diagrammed positions, we refrain from giving sketches and formal notations of matching sets of these configurations.

Also the BiTriangleX Configuration of Fig. \ref{fig:BiTriangleX}  lends itself for adding an additional group. Since in the matching set either node $d$ or $g$ need not be used, adding an additional group through $d$ and $g$ with one additional edge marker  leads to a BiTriangleX/Line Configuration, see Fig. \ref{fig:BiTriangleLines} (c).  It needs 8 markers to cover 6 groups,  a reduction of 4 nodes. Again we refrain from giving a sketch and a formal notation of a matching set of this configuration.

\begin{figure}[!h]
\begin{center}
\begin{tikzpicture}[scale=0.7]
\draw[fill=black] (-1.0,-1.732) circle (2.0pt);
\node at (-1.3,-1.732) {$a$};
\draw[fill=black] (-1.0,0.0) circle (2.0pt);
\node at (-1.3,0.0) {$b$};
\draw[fill=black] (-1.0,1.732) circle (2.0pt);
\node at (-1.3,1.732) {$c$};
\draw[fill=black] (0.5,-0.866) circle (2.0pt);
\node at (0.5,-1.166) {$d$};
\draw[fill=black] (0.5,0.866) circle (2.0pt);
\node at (0.5,1.166) {$e$};
\draw[fill=black] (2.0,-1.732) circle (2.0pt);
\node at (2.3,-1.732) {$f$};
\draw[fill=black] (2.0,0.0) circle (2.0pt);
\node at (2.3,0.0) {$g$};
\draw[fill=black] (2.0,1.732) circle (2.0pt);
\node at (2.3,1.732) {$h$};
\draw[dashed] (2.0,0.0) -- (-1.0,1.732);
\draw[dashed] (2.0,0.0) -- (-1.0,-1.732);
\draw[dashed] (-1.0,1.732) -- (-1.0,-1.732);
\draw[dashed] (-1.0,0.0) -- (2.0,-1.732);
\draw[dashed] (-1.0,0.0) -- (2.0,1.732);
\draw[dashed] (2.0,-1.732) -- (2.0,1.732);
\draw[dashed] (-1.0,0.0) -- (2.0,0.0);
\end{tikzpicture} 
\hspace*{1.0cm}
\begin{tikzpicture}[scale=0.9]
\draw[fill=black] (-1.05,-0.35) circle (2.0pt);
\draw[fill=black] (-0.35,-1.05) circle (2.0pt);
\draw[fill=black] (-0.35,-0.35) circle (2.0pt);
\draw[fill=black] (-0.35,0.35) circle (2.0pt);
\draw[fill=black] (0.35,-0.35) circle (2.0pt);
\draw[fill=black] (0.35,0.35) circle (2.0pt);
\draw[fill=black] (0.35,1.05) circle (2.0pt);
\draw[fill=black] (1.05,0.35) circle (2.0pt);
\draw[dashed] (-1.05,-0.35) -- (1.05,-0.35);
\draw[dashed] (-1.75,-1.05) -- (0.35,1.05);
\draw[dashed] (0.35,-1.05) -- (0.35,1.05);
\draw[dashed] (-0.35,-1.05) -- (-0.35,1.05);
\draw[dashed] (-1.05,0.35) -- (1.05,0.35);
\draw[dashed] (-0.35,-1.05) -- (1.75,1.05);
\draw[dashed] (1.05,-1.05) -- (-1.05,1.05);
\input{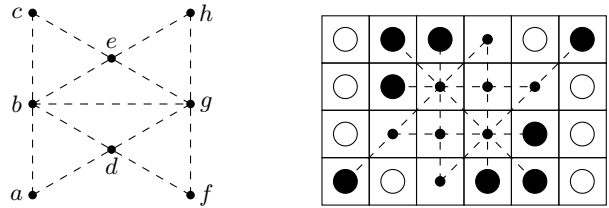}
\matrix[square matrix=.63cm]
{
	&	&	& 	& 	& 	\\
	&	&	& 	& 	& 	\\
	&	&	& 	& 	& 	\\
	&	&	& 	& 	& 	\\
};
\draw[fill=black] (-1.75,-1.05) circle (5.0pt);
\draw (-1.05,-1.05) circle (5.0pt);
\draw[fill=black] (0.35,-1.05) circle (5.0pt);
\draw[fill=black] (1.05,-1.05) circle (5.0pt);
\draw (1.75,-1.05) circle (5.0pt);
\draw (-1.75,-0.35) circle (5.0pt);
\draw[fill=black] (1.05,-0.35) circle (5.0pt);
\draw (1.75,-0.35) circle (5.0pt);
\draw (-1.75,0.35) circle (5.0pt);
\draw[fill=black] (-1.05,0.35) circle (5.0pt);
\draw (1.75,0.35) circle (5.0pt);
\draw (-1.75,1.05) circle (5.0pt);
\draw[fill=black] (-1.05,1.05) circle (5.0pt);
\draw[fill=black] (-0.35,1.05) circle (5.0pt);
\draw (1.05,1.05) circle (5.0pt);
\draw[fill=black] (1.75,1.05) circle (5.0pt);
\end{tikzpicture}
\end{center}
\vspace*{-0.2cm}
\caption{Sketch  and example position   illustrating the FlatStar/Line Configuration.}
    \label{fig:FlatStarLine}
\end{figure}

A last profitable line configuration  is the FlatStar/Line Configuration, were we have connected the two main markers $b$ and $g$ in Fig. \ref{fig:FlatStar} (left) with an additional group.  Since White always plays $b$ or $g$, we even do not need an additional marker for this group. 
Even though this Line combination is also evident, it is so efficient (the most efficient one encountered in our research so far), that we still give a sketch and an example position (see Fig. \ref{fig:FlatStarLine}).

\vspace*{-0.3cm}

\noindent
The formal notation of a matching set of this configuration is:
\begin{align*}
M = (\{&a \mhyphen h\}, \{abc,adg,bdf,beh,bg,ceg,fgh\}, \\ 
\{&a \rightarrow b \{(c,e)(d,g)(f,h)\}, 
     b \rightarrow g \{(a,c)(d,f)(e,h)\}, \\
   &d \rightarrow b \{(a,g)(c,e)(f,h)\}\}) 
\end{align*}

\noindent were we again have omitted the equivalent coverings when Black starts with $c$, $e$, $f$, $g$, or $h$.
This configuration is consequently even more efficient, needing just 8 markers to cover 7 groups, a reduction of 6 nodes.

\subsection{PolyCycle Configurations  with One Common Side \label{subsec:PolyCommon}}

\begin{figure}[!h]
\begin{center}
\begin{tikzpicture}[scale=0.7]
\draw[fill=black] (-1.0,-1.732) circle (2.0pt);
\node at (-1.0,-2.032) {$a$};
\draw[fill=black] (-1.0,1.732) circle (2.0pt);
\node at (-1.0,2.032) {$b$};
\draw[fill=black] (2.0,-1.732) circle (2.0pt);
\node at (2.3,-2.032) {$c$};
\draw[fill=black] (0.5,-1.732) circle (2.0pt);
\node at (0.5,-2.032) {$d$};
\draw[fill=black] (0.5,0.0) circle (2.0pt);
\node at (0.8,0.0) {$e$};
\draw[fill=black] (2.0,1.732) circle (2.0pt);
\node at (2.3,2.032) {$f$};
\draw[fill=black] (0.5,1.732) circle (2.0pt);
\node at (0.5,2.032) {$g$};
\draw[fill=black] (-4.0,0.0) circle (2.0pt);
\node at (-4.3,0.0) {$h$};
\draw[fill=black] (-2.5,-0.866) circle (2.0pt);
\node at (-2.65,-1.166) {$i$};
\draw[fill=black] (-2.5,0.866) circle (2.0pt);
\node at (-2.65,1.166) {$j$};
\draw[dashed] (2.0,-1.732) -- (-1.0,-1.732); 
\draw[dashed] (2.0,-1.732) -- (-1.0,1.732); 
\draw[dashed, ultra thick] (-1.0,1.732) -- (-1.0,-1.732); 
\draw[dashed] (2.0,1.732) -- (-1.0,-1.732); 
\draw[dashed] (2.0,1.732) -- (-1.0,1.732); 
\draw[dashed] (-4.0,0.0) -- (-1.0,-1.732); 
\draw[dashed] (-4.0,0.0) -- (-1.0,1.732); 
\end{tikzpicture} 
\hspace*{1.0cm}
\begin{tikzpicture}[scale=0.9]
\draw[fill=black] (0.0,-0.7) circle (2.0pt); 
\draw[fill=black] (0.0,0.7) circle (2.0pt); 
\draw[fill=black] (1.4,0.7) circle (2.0pt); 
\draw[fill=black] (0.7,0.0) circle (2.0pt); 
\draw[fill=black] (0.7,0.7) circle (2.0pt); 
\draw[fill=black] (0.7,-0.7) circle (2.0pt); 
\draw[fill=black] (1.4,-0.7) circle (2.0pt); 
\draw[fill=black] (-1.4,-0.7) circle (2.0pt); 
\draw[fill=black] (-1.4,0.7) circle (2.0pt); 
\draw[fill=black] (-0.7,0.0) circle (2.0pt); 
\draw[dashed, ultra thick] (0.0,-0.7) -- (0.0,1.4);
\draw[dashed] (0.0,0.7) -- (2.1,0.7);
\draw[dashed] (-0.7,-1.4) -- (1.4,0.7);
\draw[dashed] (0.0,-0.7) -- (2.1,-0.7);
\draw[dashed] (1.4,-0.7) -- (-0.7,1.4);
\draw[dashed] (-2.1,-1.4) -- (0.0,0.7);
\draw[dashed] (0.0,-0.7) -- (-2.1,1.4);
\input{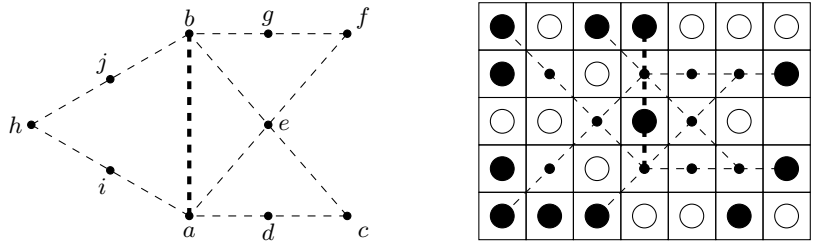}
\matrix[square matrix=.63cm]
{
	&	&	&	& 	& 	& 	\\
	&	&	&	& 	& 	& 	\\
	&	&	&	& 	& 	& 	\\
	&	&	&	& 	& 	& 	\\
	&	&	&	& 	& 	& 	\\
};
\draw[fill=black]  (-2.1,-1.4) circle (5.0pt); 
\draw[fill=black]  (-1.4,-1.4) circle (5.0pt); 
\draw[fill=black]  (-0.7,-1.4) circle (5.0pt); 
\draw  (0.0,-1.4) circle (5.0pt); 
\draw  (0.7,-1.4) circle (5.0pt); 
\draw[fill=black]  (1.4,-1.4) circle (5.0pt); 
\draw (2.1,-1.4) circle (5.0pt); 
\draw[fill=black] (-2.1,-0.7) circle (5.0pt); 
\draw  (-0.7,-0.7) circle (5.0pt); 
\draw[fill=black] (2.1,-0.7) circle (5.0pt); 
\draw (-2.1,0.0) circle (5.0pt); 
\draw  (-1.4,0.0) circle (5.0pt); 
\draw[fill=black] (0.0,0.0) circle (5.0pt); 
\draw  (1.4,0.0) circle (5.0pt); 
\draw[fill=black]  (-2.1,0.7) circle (5.0pt); 
\draw  (-0.7,0.7) circle (5.0pt); 
\draw[fill=black]  (2.1,0.7) circle (5.0pt); 
\draw[fill=black]  (-2.1,1.4) circle (5.0pt); 
\draw  (-1.4,1.4) circle (5.0pt); 
\draw[fill=black]  (-0.7,1.4) circle (5.0pt); 
\draw[fill=black]  (0.0,1.4) circle (5.0pt); 
\draw  (0.7,1.4) circle (5.0pt); 
\draw  (1.4,1.4) circle (5.0pt); 
\draw  (2.1,1.4) circle (5.0pt); 
\end{tikzpicture}
\hspace*{1.0cm}
\end{center}
\vspace*{-0.2cm}
\caption{Sketch  and example position  illustrating the TriTriangleX Configuration.}
    \label{fig:TriTriangleX}
\end{figure}

It is possible to include more than two cycles in a configuration with all cycles having one common side. 
We give just one example, namely a TriTriangleX Configuration sketched in Fig. \ref{fig:TriTriangleX}. 

\vspace*{-0.3cm}

\noindent
The formal notation for a matching set of this configuration is:
\begin{align*}
M = (\{&a \mhyphen j\}, \{ ab,adc,aef,aih,bec,bgf,bjh\}, \\
\{ &a \rightarrow b\{(c,d),(e,f),(h,i)\}, 
      c \rightarrow a\{(b,e),(f,g),(h,j)\}, \\
   &d \rightarrow a\{(b,f),(c,e),(h,j)\}, 
     e \rightarrow a\{(b,c),(f,g),(h,j)\}, \\
   &h \rightarrow a\{(b,j),(c,e),(f,g)\}, 
      i \rightarrow a\{(b,c),(f,g),(h,j)\}\}) 
\end{align*}

\noindent were we have omitted the equivalent coverings when Black starts with  $b$, $f$, $g$, or $j$. 
In this configuration the 10 markers cover 7 groups, a reduction of 4 nodes. 
An example position is given in the  diagram in Fig. \ref{fig:TriTriangleX}. 


\section{Overview of Reductions for Configurations}

In  Table \ref{tab:NewConfigs} we show the requirements and reductions for all configurations given in Sections 3-5. Moreover, we give the results of a small experiment on the example positions.
\noindent The first column gives the configuration, the next two columns the number of markers required and the number of groups. Then the reductions in nodes required compared to HJ-pairings are given and the coverage ratios (markers/groups). In the last two columns we give the example positions used and the number of nodes that our 4-in-a-Row program \cite{Uiterwijk2017} needs to prove them to be draws, using naive HJ-pairings. 

\begin{table}[!h]
\centering
\begin{tabular}{l | c | c | c | c || c | c }
\hline
Configuration 	& $\#$ markers	& $\#$ groups	& $\#$ reduction	&coverage ratio	& diagram				& $\#$ nodes \\	
\hline
\hline
\multicolumn{5}{l}{Single cycles} 										\\
\hline
Triangle		& 5			& 3			& 1			& 1.67		& Fig. \ref{fig:Triangle}		& 17	\\
Square		& 7			& 4			& 1			& 1.75		& Fig. \ref{fig:Square}		& 40	\\
Cycle $C_n$		& $2n-1$		& $n$			& 1			& $(2n-1) / n$	& - 					& - 	\\
\hline
\hline
\multicolumn{5}{l}{BiCycles with 1 new side} 									\\
\hline
Triangle/Line		& 6			& 4			& 2			& 1.50		& Fig. \ref{fig:TriangleLine}	& 26	\\
Square/Line		& 8			& 5			& 2			& 1.60		& Fig. \ref{fig:SquareLine}	& 67	\\
Cycle $C_n$/Line	& $2n$		& $n+1$		& 2			& $2n / (n+1)$	& - 					& - 	\\
\hline
\hline
\multicolumn{5}{l}{BiCycles with 2 new sides} 								\\
\hline
BiTriangle		& 8			& 5			& 2			& 1.60		& Fig. \ref{fig:BiTriangle}		& 25	\\
BiTriangleX		& 7			& 5			& 3			& 1.40		& Fig. \ref{fig:BiTriangleX}		& 17	\\
\hline
\hline
\multicolumn{5}{l}{BiCycles with 3 new sides}									\\
\hline
FlatStar		& 8			& 6			& 4			& 1.33		& Fig. \ref{fig:FlatStar}		& 17	\\
\hline
\hline
\multicolumn{5}{l}{PolyCycles} 						\\
\hline
BiTriangle/Line	& 9			& 6			& 3			& 1.50		& Fig. \ref{fig:BiTriangleLines}a	& 25	\\
BiTriangle/BiLine	& 10			& 7			& 4			& 1.43		& Fig. \ref{fig:BiTriangleLines}b	& 25	\\
BiTriangleX/Line	& 8			& 6			& 4			& 1.33		&Fig. \ref{fig:BiTriangleLines}c	& 51	\\
FlatStar/Line		& 8			& 7			& 6			& 1.14		& Fig. \ref{fig:FlatStarLine}	& 17	\\
TriTriangleX		& 10			& 7			& 4			& 1.43		& Fig. \ref{fig:TriTriangleX}		& 23	\\
\hline
\end{tabular}
\vspace*{0.2cm}
\caption{Some drawing configurations.  
\label{tab:NewConfigs}}
\end{table}


\section{Conclusions and Future Research}
In this paper we have presented an enhancement of the Hales-Jewett pairing strategy, called {\em Set Matching}. This strategy needs  less than 2 squares per group. 
We have given several configurations with their matching sets, including Cycle, BiCycle, and PolyCycle Configurations. Depending on the configuration, the coverage ratio is reduced to 1.14, compared to 2.0 for HJ-pairings.
Many example positions for $4$-in-a-Row games are given, including the direct proof (without further search) that the $4 \times 4$ board is a draw.

As future research  we first will implement Set Matching in our $k$-in-a-Row solver \cite{Uiterwijk2017} to investigate its effectiveness compared to naive HJ-pairings for solving games. 
Secondly, we will further look for other interesting  configurations, a main challenge being one with a matching set ($N$, $G$, $C$) with $|N| \le |G|$, i.e., a coverage factor of 1.0 or lower.  
Thirdly, we will investigate how well Set Matching will perform in case of overlapping groups. 



\begin{thebibliography}{[MT1]}

\bibitem{Beck2005}
Beck, J.: 
Positional Games.
{\it Combinatorics, Probability and Computing} {\bf 14} (2005), 649--696. 

\bibitem{Beck2008}
Beck, J.: 
{\it Combinatorial Games: Tic-Tac-Toe Theory.} 
Cambridge University Press, 2008.

\bibitem{Gardner1983}
Gardner, M.:
Ticktacktoe Games. 
Ch. 9 in {\em Wheels, Life, and Other Mathematical Amusements}. 
W.H. Freeman, 1983, 94--105. 

\bibitem{HalJew1963}
Hales, A.W., and Jewett, R.I.:
Regularity and Positional Games.
{\it Trans. Amer. Math. Soc.} {\bf 106} (1963) 222--229.

\bibitem{Uiterwijk2014}
Uiterwijk, J.W.H.M.:
Perfectly Solving {D}omineering Boards.
{\em Computer Games. Workshop on Computer Games, CGW 2013, Beijing, China. Revised Selected Papers}
(eds. T. Cazenave, M.H.M Winands, and H. Iida),
{\em Communications in Computer and Information Science} {\bf 408}  (2014) 97--121, Springer.

\bibitem{Uiterwijk2017}
Uiterwijk, J.W.H.M.:
Solving Four-in-a-Row.
Submitted (2017).

\bibitem{UitHer2000}
Uiterwijk, J.W.H.M. and Herik, H.J. van den:
The Advantage of the Initiative.
{\it Information Sciences} {\bf 122}(1) (2000) 43--58.
\end{thebibliography}
\end{document}